\documentclass[11pt]{amsart}
\usepackage{amsmath, tikz}
\usepackage{color}

\newcommand{\dx}{~\mathrm{d}x}
\newcommand{\sym}{\mathrm{sym}}

\newcommand{\Id}{\mathrm{Id}}

\newcommand{\R}{\mathbb{R}}

\def\be{\begin{equation}}
\def\ee{\end{equation}}

\begin{document}
\title[Isometric immersions]{Isometric immersions and applications}
\author{Qing Han}
\address{Q.H.: Department of Mathematics, University of Notre Dame, Notre Dame, IN 46556}
\author{Marta Lewicka}
\address{M.L.: University of Pittsburgh, Department of Mathematics, 
139 University Place, Pittsburgh, PA 15260}
\email{Qing.Han.7@nd.edu, lewicka@pitt.edu} 

\thanks{Q.H. was supported by NSF grant DMS-2305038. M.L. was
  supported by NSF grant DMS-2006439. The authors thank Marcus
Khuri for the assistance in writing section 4.} 

\begin{abstract}
We provide an introduction to the old-standing problem of isometric
immersions. We combine a historical account of its multifaceted advances, which
have fascinated geometers and analysts alike, with some of the
applications in the mathematical physics and mathematical materials science,
old and new.
\end{abstract}

\maketitle

\section{The isometric immersion problem}

\noindent The concept of a {\em Riemannian manifold} $(M^d, g)$, an abstract
$d$-dimensional manifold with a metric structure, was first formulated by Bernhard Riemann in 1868 to
generalize the classical objects such as curves and surfaces in $\R^3$. 
A manifold is a topological space that locally resembles Euclidean space near each point:
in $M$ (we write $M^d$ only to emphasize the dimension) each point
$p\in M$ has a neighborhood that is homeomorphic to an open subset of 
$\mathbb R^d$. Induced by that homeomorphism, there is a tangent space
$T_pM$, namely a $d$-dimensional vector space
that is further equipped with a positive-definite inner product $g_p$. 
The family $g=\{g_p\}_{p\in M}$ is called a {\em Riemannian metric} on $M$.
For details, we recommend consulting the textbook \cite{HH} and
references listed therein.

\medskip

\noindent All Euclidean spaces are manifolds and, endowed with the standard
Euclidean inner product, are Riemannian manifolds.  
The simplest nontrivial $d=2$ - dimensional example is the unit sphere $\mathbb S^2$ 
consisting of all unit vectors in $\mathbb R^3$. 
At each $p\in \mathbb S^2$, the tangent space $T_p\mathbb S^2$ is the 
hyperplane in $\mathbb R^3$ perpendicular to the unit vector $p$. 
The standard inner product in $\mathbb R^3$ restricted to $T_p\mathbb
S^2$ induces then a Riemannian metric on $\mathbb S^2$, called the
round metric. Naturally, there arises the
question of whether any abstract $(M^d,g)$ can be identified, in a
similar manner, as a submanifold of some Euclidean space $\R^n$ with its induced
metric. This is the {\em isometric embedding question}, which has assumed a position of
fundamental conceptual importance in Differential Geometry. 

\medskip

\noindent The metric requirement can be expressed in terms of partial
differential equations (PDEs). Consider a special case:
let $B_1$ be the unit ball in $\mathbb R^d$, regarded as a
$d$-dimensional manifold. The given inner product $g_p$ can be
identified with a $d\times d$ symmetric, positive definite matrix at each point $p\in
B_1$, so that $g$ is a function from $B_1$ to $\R^{d\times d}_{\sym, >}$.
Isometrically {immersing} $(B_1, g)$ in some Euclidean space $\mathbb R^n$ 
means that there exists $u=(u^1,\ldots,u^n):B_1\to \mathbb R^n$ such
that the induced metric agrees with $g$ at each point:
\begin{equation}\label{A.1.4}
(\nabla u)^T\nabla u = g \quad\text{in }B_1.
\end{equation}
Above, we used the matrix notation for the transpose and product. 
We call the map  $u$ in (\ref{A.1.4}) an {\em isometric embedding} or {\em immersion}
according to whether it is injective or not. 
In the general setting, and when globally defined,  both sides of
\eqref{A.1.4} become second order covariant tensors. 
We now briefly review several important aspects of isometric
embeddings and immersions of Riemannian manifolds in the Euclidean space.

\section{General isometric embeddings and immersions of Riemannian manifolds}

\subsection{The analytic solutions.} 
Let us examine \eqref{A.1.4} closely. First, we notice
that it has $n$ unknowns $\{u^i\}_{i=1}^n$ and there are $s_d\doteq
{d(d+1)}/2$ equations, corresponding to the entries in $d\times d$ symmetric matrices. 
The integer $s_d$ is called the {\em Janet dimension}.
As a rule in solving PDEs, the number of unknowns 
should be bigger than or equal to the number of equations,
otherwise, solutions are not expected to exist in general. Hence we
require that $n\ge s_d$. In 1873, Ludwig Schlaefli conjectured that every
$d$-dimensional smooth Riemannian manifold admits a {\em smooth
  local isometric embedding} in $\mathbb R^{s_d}$. It
was more than 50 years later that an affirmative answer was given
in the analytic case by Maurice Janet (for $d=2$) and
\'Elie Cartan (for $d\ge 3$). 

\subsection{The smooth case and the Nash-Moser iteration.} 
Second, \eqref{A.1.4} is classified as a first-order nonlinear 
system of PDEs, the order of derivatives being $1$ and the derivatives
of $u$ appearing quadratically (nonlinearly) in the system. 
In general, nonlinear PDEs are solved by Newton's method, 
an iteration process generating an actual solution from an
approximate one. Here, the difficulty arises from the loss of
derivatives at each iteration. 

\medskip

\noindent To explain this clearly, recall how Newton's method is used
to find a root of a single-variable function.  
Let $f$ be a scalar-valued twice continuously differentiable function on $\mathbb R$. 
Pick an initial guess $t_0\in\mathbb R$ for a root of $f$, meaning that $f(t_0)$ is small. 
To find $t_1=t_0+s_1$ such that $f(t_1)$ is smaller, we use Taylor's theorem to have:
$$f(t_1)=f(t_0)+f'(t_0)s_1+\frac12f''(t_0+\xi s_1)s_1^2
\quad\text{for some } \xi\in (0,1).$$
It is natural to choose $s_1$ such that $f(t_0)+f'(t_0)s_1=0$ and hence:
$$s_1=-\frac{f(t_0)}{f'(t_0)}, \qquad f(t_1)=\frac12f''(t_0+\xi s_1)s_1^2.$$
We see that $s_1=t_1-t_0$ is linear in $f(t_0)$ and so $f(t_1)$ is quadratic in $f(t_0)$. 
We can iterate this process and obtain a sequence $\{t_l\}_{l\to\infty}$ inductively by:
\begin{equation}\label{eq-choice-s-l}
\begin{split}
t_l\doteq t_{l-1}+s_l, \qquad s_l\doteq -\frac{f(t_{l-1})}{f'(t_{l-1})},\qquad 
f(t_l)=\frac12f''(t_{l-1}+\xi_l s_l)s_l^2.
\end{split}
\end{equation}
We next need to prove that $\{t_l\}_{l\to\infty}$ converges to some $t_\infty$ and
that $\{f(t_l)\}_{l\to\infty}$ converges to $0$. As a consequence, $f(t_\infty)=0$ resulting in a root of $f$
at $t_\infty$. An important condition in the convergence proof is given by:
\begin{equation*}\label{eq-condition-derivative-l}
|f'(t_l)|\ge c\quad\text{for all } l\geq 0,
\end{equation*}
for some positive constant $c$, independent of $l$. 
With \eqref{eq-condition-derivative-l} and an appropriate condition on $f''$, 
we get that indeed $\{s_l\}_{l\to\infty}$ converges to zero fast and so
$f(t_l)\to 0$.

\medskip

\noindent We now apply the same idea to \eqref{A.1.4}, equivalent to finding roots of:
$$G[u]=(\nabla u)^T\nabla u - g.$$
We first take a $u_0: B_1\to \mathbb R^n$ such that 
$G[u_0]$ is small, write $u_1=u_0+v_1$ and compute:
$$G[u_1]=G[u_0]+ 2\,\sym \big((\nabla u_0)^T\nabla v_1\big) + (\nabla v_1)^T\nabla v_1,$$
where we arranged the expression in the right-hand side according to
the powers of (the derivatives of) $v_1$. As in the scalar case, we
would like to choose $v_1$ such that:
\begin{align}\label{eq-linear-system-0}
G[u_0] + (\nabla u_0)^T\nabla v_1 + (\nabla v_1)^T\nabla u_0=0  \quad\text{in } B_1.
\end{align}
This is a first-order linear system of PDEs that can
be solved by imposing appropriate boundary values. 
However, there is a serious issue in the iteration process. 

\medskip

\noindent In studying PDEs or systems thereof,
we need to introduce appropriate spaces of functions and equip them with appropriate norms 
to form Banach spaces. The process outlined above 
is performed in Banach spaces to generate a \lq\lq root\rq\rq\ of the functional $G$. 
An important feature of \eqref{eq-linear-system-0} is that
although a linear combination of derivatives of $v_1$ is governed by $G[u_0]$,  
it is $v_1$, not $\nabla v_1$, that is determined by $G[u_0]$. 
More generally, derivatives of $G[u_0]$ of a certain order determine only derivatives of $v_1$ 
up to the same order. Hence, in the inductive version of
\eqref{eq-linear-system-0} and a counterpart of \eqref{eq-choice-s-l}:
\begin{equation}\label{eq-linear-system-l}
\begin{split}
u_l\doteq u_{l-1}+v_l, \qquad G[u_{l-1}]+ 2\,\sym\big((\nabla
u_{l-1})^T\nabla v_l\big) =0, \qquad G[u_l]=(\nabla v_l)^T\nabla v_l,
\end{split}
\end{equation}
each step in the iteration contributes 
a loss of derivatives. If $u_0$ is, say, $100$ times continuously
differentiable then $G[u_0]$, as it involves $\nabla u_0$, is 
$99$ times continuously differentiable and thus
$v_1$ is also $99$ times differentiable and so is $u_1$. 
Then, $u_{100}$ will be continuous only and $G[u_{100}]$ does not make sense. 
As a consequence, the iteration process is terminated. 

\medskip

\noindent In an outstanding paper
published in 1956, John Nash introduced an important technique of {\em
  smoothing operators} to compensate for the aforementioned loss of derivatives. He
proved that any smooth $d$-dimensional Riemannian manifold admits
a (global) smooth isometric embedding  in 
$\mathbb R^n$, for $n=3s_d+4d$ in the compact case and
$n=(d+1)(3s_d+4d)$ in the general case. 
Specifically for \eqref{eq-linear-system-l}, Nash replaced $u_{l-1}$ by 
a function with a better regularity to regain the lost derivative,
solved for $v_l$, and modified the iteration process accordingly. 
His technique proved to be extremely useful in nonlinear
PDEs, and it is now known as
 the {\em hard implicit function theorem}, or the {\em Nash-Moser iteration}.

\subsection{Quest for the smallest dimension \textit{\textbf{n}}.}  

Following Nash, one naturally looks for the smallest
dimension of the ambient space. In 1970, Mikhael Gromov and Vladimir
Rokhlin, and independently John Greene, proved 
that any $d$-dimensional smooth Riemannian manifold admits a local smooth
isometric embedding in $\mathbb{R}^{s_d+d}$. The
proof is based on Nash's iteration scheme. In his book \cite{Gromov},  Gromov
studied various problems related to the isometric embedding of
Riemannian manifolds. He proved that $n=s_d+2d+3$ is enough for the
compact case. Then in 1989, Matthias G\"{u}nther vastly simplified Nash's
original proof: by rewriting the differential equations cleverly, he
was able to employ the contraction mapping principle, instead of the
Nash-Moser iteration, to construct solutions. G\"{u}nther also
improved the dimension of the target space to $n=\max
\{s_d+2d,s_d+d+5 \}$.  It is still not clear whether this is the best possible result.

\medskip

\noindent For $d=2$ better results are
available. According to Gromov and G\"{u}nther, any compact
$2$-dimensional smooth Riemannian manifold can be isometrically
embedded in $\mathbb R^{10}$ smoothly, whereas the local version in
$\R^4$ is due to Eduard Poznyak in 1973.  The case of $d=2$ and $n=s_2=3$ will be
discussed in a separate section below.

\medskip

\noindent On the other hand, the case $d\ge 3$, $n=s_d$ is
sharply different from the $2$-dimensional case in which there is only one
curvature function. This function determines the type of the Darboux equation
associated to (\ref{A.1.4}), whereas for $d\ge 3$,
the role of various curvatures is not clear.
In 1983, Robert Bryant, Phillip Griffiths and Dean Yang studied the characteristic
varieties associated with differential systems for the
isometric embedding in $\mathbb R^{s_d}$ of smooth $(M^d,g)$. They proved that these 
varieties are never empty if $d\ge 3$, implying, in particular,
that the governing systems are never elliptic, no matter what assumptions are put on
curvatures. They also proved that the characteristic varieties are smooth for $d =
4$ and not smooth for $d = 6, 10, 14, \ldots$. 
In 2012, Qing Han and Marcus Khuri extended that result to any $d\ge 5$ under an
additional \lq\lq smallness\rq\rq\, assumption.  

\medskip

\noindent For $d=3$, Bryant, Griffiths and Yang in 1983 classified the
type of differential systems for the isometric embedding by its
curvature functions. Here, an important quantity is the signature of
the curvature tensor viewed as a symmetric linear operator acting on
the space of $2$-forms. In particular, any smooth $3$-dimensional
Riemannian manifold admits a smooth local isometric embedding in
$\R^6$ if the signature is different from $(0,0)$ and $(0,1)$.
In 1989, Yusuke Nakamura and Yota Maeda showed the same existence
result if the curvature tensors are not zero; their key argument 
was the local existence of solutions to
differential systems of principal type. In 2018, Chen, Clelland,
Slemrod, Wang and Yang provided an alternative proof using strongly symmetric positive systems.

\subsection{Non-smooth immersions.}

On the other end of the spectrum, there are results showing that
isometric immersions have completely different qualitative behaviours
at low and high regularity. Nash in 1954 and Nicolaas Kuiper in
1955 proved the existence of a global $C^1$ isometric embedding of
$d$-dimensional Riemannian manifolds in $\mathbb R^{d+1}$. These are
not mere existence statements, as their results in fact show that 
every {\em short immersion} (or embedding),
i.e. $u: B_1\to\R^{d+1}$ for which condition (\ref{A.1.4}) is replaced by:
\begin{equation}\label{B.1}
(\nabla u)^T\nabla u < g \quad\text{in }B_1,
\end{equation}
can be uniformly approximated by $C^1$-regular actual
solutions (immersions or embeddings) to (\ref{A.1.4}). 
The inequality above is understood pointwise, in the sense of 
matrices. This abundance of solutions, usually referred to as {\em
  flexibility results}, is typical in applications of Gromov's {\em
  h-principle} in which a PDE is replaced by a partial differential
relation (a differential inclusion) whose solutions are then modified
through an iteration technique called {\em convex integration} to
produce a nearby solution of the underlying PDE.

\medskip

\noindent In 1965, Yuri Borisov used this approach to H\"older - regular solutions and announced
that flexibility holds with regularity $C^{1,\alpha}$ for
any $\alpha<\frac{1}{1+2 s_d}$ and analytic $g$ on $B_1$. He
subsequently gave full details of the proof in 2004 for
dimension $d=2$ and $\alpha<\frac{1}{13}$. In 2012 Sergio Conti,
Camillo De Lellis and Laszlo Szekelyhidi validated the original
Borisov's statement in case of $C^2$ metrics on $d$-dimensional balls,
and in case of compact Riemannian manifolds $(M^d,g)$ with $\alpha<\frac{1}{1+2s_d(d+1)}$.
The same results and exponent bounds also hold for any target ambient
dimension $n>d$. However, in 1978
Anders K\"allen proved that any $C^\beta$ metric, with $\beta<2$,
allows for flexibility up to exponent $\frac{\beta}{2}$, provided that
$n$ is sufficiently large. Recently, the second author of this paper
proved flexibility of the related {\em Monge-Ampere system}, which is
the linearization of the isometric immersion
problem (\ref{A.1.4}) around $g=\Id_d$, and in which any $C^1$ -
regular subsolution can be uniformly approximated by $C^{1,\alpha}$
exact solutions, for any $\alpha<\frac{1}{1+2s_d/(n-d)}$, in agreement
with the fully nonlinear case at $n=d+1$ and the Kallen result when
$n\to\infty$. 

\medskip

\noindent Dependence of the flexibility threshold exponents on
$s_d$ reflects the technical limitation of the method rather than the absolute lack of flexibility beyond those
thresholds. In the proofs, the symmetric, positive definite ``defect''
$\mathcal{D}=g-(\nabla u)^T\nabla u$ is decomposed  
into a linear combination of $s_d$ rank-one defects with
nonnegative coefficients. Each of these ``primitive defects'' is then
cancelled by adding to $u$ a small but fast oscillating perturbation,
which however causes increase of the second derivative of $u$ by the factor
of the oscillation frequency. The ultimate H\"older regularity of the
approximating immersion interpolates between the controlled $C^1$
norms of the adjusted $u$ and the blow-up rate of the $C^2$ norm, dictated by the
number of these one-dimensional adjustments. In case of higher
codimension when $n>d+1$, several
primitive defects may be cancelled at once, reducing the
blow-up rate of second derivatives and thus improving the regularity exponent $\alpha$.
In the same vein, if the number of
primitive defects in the decomposition of $\mathcal{D}$ could be
lowered, for example by an appropriate change of 
variables, then flexibility would hold with higher $\alpha$. This
observation is precisely behind the improved regularity statements for
$d=2$ - dimensional problems, listed in the next section.

\section{Isometric embeddings of surfaces}

\subsection{Local isometric embedding of surfaces in $\mathbb R^3$.}
We now give an overview of the 
question of isometrically embedding a $2$-dimensional
Riemannian manifold in $\mathbb R^3$.
There are basically two methods to study the local case.
The first one, already known to Jean Darboux in 1894, restates the
problem equivalently as finding a local
solution of a nonlinear equation of the {\em Monge-Amp\`{e}re type}. 
Specifically, let $g$ be a $C^r$-metric on a simply
connected $\Omega\subset \mathbb R^2$ for some $r\in
[2,\infty]$. If there exists $u\in C^s(\omega,\R)$, solving:
\begin{equation}\label{1.2.3}
\det(\nabla_g^2{u})=K(\det g)(1-|\nabla_g u|^2),
\end{equation}
with  $|\nabla_gu|<1$ for some $s\in [2, r]$, then $(\Omega, g)$
admits a $C^s$-isometric immersion in $\mathbb R^3$.
The equation \eqref{1.2.3} is now called the {\em Darboux equation}; and its type is
determined by the sign of the Gauss curvature $K$ of $g$: elliptic if
$K$ is positive, hyperbolic if $K$ is negative, and  degenerate
if $K$ vanishes. Remarkably, even today, the local solvability of
the Darboux equation in the general case is not covered by any known
theory of PDEs.

\medskip

\noindent A different method to study the local isometric embedding of
surfaces in $\mathbb{R}^3$ relies on
the classical theory of surfaces asserting that such
immersion exists provided the solvability of the {\em Gauss-Codazzi
  system}. Namely, let $\{\Gamma^i_{jk}\}_{i,j,k=1,2}$ be the
Christoffel symbols of the given metric $g$,  and $K$ its Gauss
curvature. Then the coefficients of the second fundamental form
$I\!\!I= Ldx_1^2+2Mdx_1dx_2+Ndx_2^2$ satisfy:
\begin{equation}\label{e:GC}
\begin{split}
& \partial_{2}L-\partial_{1}M=\Gamma^1_{12}L+(\Gamma^2_{12}-\Gamma^1_{11})M-\Gamma^2_{11}N,\\
&\partial_{2}M-\partial_{1}N=\Gamma^1_{22}L+(\Gamma^2_{22}-\Gamma^1_{21})M-\Gamma^2_{21}N,\\
& LN-M^2=K(g_{11}g_{22}-g_{12}^2).
\end{split}
\end{equation}
We note in passing that the first attempt to establish the local isometric embedding of
surfaces in $\mathbb R^3$ was neither through (\ref{1.2.3}) nor
(\ref{e:GC}): in 1908, Hans Levi solved the case of 
surfaces with negative curvature by using the equations of virtual
asymptotes. 

\medskip

\noindent It was several decades later that (\ref{1.2.3})
attracted attention of those interested in the isometric
embedding. In the early 1950s, Philip Hartman and Aurel Wintner
studied (\ref{1.2.3}) with $K\neq 0$ and proved existence of its local solution.
The case when $K$ vanishes did not give way to the
efforts of mathematicians for a long time. In 1985 and 1986, Chang-Shou Lin made important
breakthroughs, establishing existence in a neighborhood of
$p\in\Omega$ such that  $K(p)=0$ and $dK(p)\ne 0$ (in 2005 Han gave an alternative proof of
this result), or $p$ such that $K\ge 0$ in the whole neighbourhood. Later, in 1987,
Gen Nakamura covered the case of $K(p)=0$, $dK(p)=0$ and 
$\mbox{Hess}\, K(p)<0$. For the case of nonpositive $K$, 
Jia-Xing Hong in 1991 also proved the existence of a sufficiently
smooth local isometric embedding in a neighborhood of $p$ if $K=h
\varphi^{2m}$, where $h$ is a negative function and $\varphi$ is a function with
$\varphi(p)=0$ and $d\varphi(p)\neq 0$.
In 2010, Han and Khuri proved existence of the smooth local isometric embedding near $p$ 
if $K$ changes its sign only along two smooth curves 
intersecting transversely at $p$. 
All these results are based on careful studies of the Darboux equation.

\medskip

\noindent In 2003, Han, Hong, and Lin studied (\ref{e:GC}) and proved the local
isometric embedding for a large class of metrics with nonpositive
Gauss curvature $K$, for which directional derivative has a simple
structure for its zero set. This gives the results of
Nakamura and Hong as special cases.
On the other hand, Aleksei Pogorelov in 1972 constructed a  $C^{2,1}$
metric $g$ on $B_{1}\subset \mathbb{R}^{2}$ with a sign-changing
$K$ such that $(B_{r}, g)$ cannot be realized as a
$C^{2}$ surface in $\mathbb{R}^{3}$ for any $r>0$.
Nikolai Nadirashvili and Yu Yuan in 2008 constructed a $C^{2,1}$
metric $g$ on $B_{1}$ of $K\geq 0$, 
with no $C^2$ local isometric embedding in $\mathbb{R}^{3}$.

\subsection{Global isometric embedding of surfaces in $\mathbb R^3$.}

In 1916, Hermann Weyl posed the following problem: does every smooth metric
on $\mathbb S^2$ with positive Gauss curvature admit a
smooth isometric embedding in $\mathbb R^3$? The first attempt to
solve the problem was made by Weyl himself, through the continuity
method and the a priori estimates up to the second derivatives.
Twenty years later, Hans Lewy solved the problem for analytic metrics
$g$. In 1953, Luis Nirenberg gave a complete solution under the very mild
hypothesis that $g$ is $C^4$. The result was extended to the $C^3$ case 
by Erhard Heinz in 1962. In a completely different approach,  
Aleksandr Alexandroff in 1942 obtained a generalized solution
to Weyl's problem as a limit of polyhedra. Further,
it is known from the work of Pogorelov in the 1950s that closed $C^1$
surfaces with positive Gauss curvature and bounded extrinsic
curvature are convex, and that closed convex surfaces are {\em rigid}
in the sense that their isometric immersions are unique up to rigid motions. 
In 1994 and 1995, Pei-Feng Guan and Yanyan Li, and Hong and Zuily independently
generalized Nirenberg's result for metrics on $\mathbb S^2$ with
nonnegative Gauss curvature.

\medskip

\noindent The investigation of the isometric immersion of metrics with negative
curvature goes back to David Hilbert. He proved in 1901
that the full hyperbolic plane cannot be isometrically immersed in
$\mathbb R^3$. The next natural step is to extend such a result to
complete surfaces whose Gauss curvature is bounded above by a
negative constant.  The final solution to this problem was obtained
by Nikolai Efimov in 1963: he proved that
any complete negatively curved smooth surface does not admit a $C^2$
isometric immersion in $\mathbb R^3$ if its Gauss curvature is
bounded away from zero. In the years following, Efimov extended his result in
several ways.
Before the 1970s, the study of negatively curved surfaces was
largely directed at the nonexistence of isometric immersions in $\mathbb
R^3$. In the 1980s,  Shing-Tung Yau proposed to find a sufficient
condition for complete negatively curved surfaces to be
isometrically immersed in $\mathbb R^3$.  In 1993,  Hong identified such
condition in terms of the Gauss curvature decaying at a certain
rate at infinity.  His discussion was based on a differential system
equivalent to the Gauss-Codazzi system (\ref{e:GC}).

\subsection{Non-smooth immersions of surfaces.} 

Extension of the rigidity of Weyl's problem to H\"older-regular $C^{1,\alpha}$ isometric
immersions, is originally due to Borisov in the late 1950s, who proved
that for $\alpha>\frac{2}{3}$, the image of a surface with positive
Gauss curvature has bounded extrinsic curvature. Hence, if
$(\mathbb{S}^2, g)$ is a Riemannian manifold with $K>0$ then
$u(\mathbb{S}^2)$ is the boundary of a bounded convex set that is
unique up to rigid motions of $\R^3$, provided that $u\in
C^{1,\alpha}$ with $\alpha>\frac{2}{3}$. In particular, if $K$ is constant
then $u(\mathbb{S}^2)=\frac{1}{\sqrt{K}}\mathbb{S}^2$. In 2012 Conti, De Lellis and Szekelyhidi
provided a direct analytic proof of these results, based on a clever use
of the commutator estimate. 

\medskip

\noindent As we have mentioned before, flexibility for isometric
immersions of surfaces in $\R^3$ has been proved by Conti, De Lellis and Szekelyhidi, up to
the regularity exponent $\frac{1}{7}$. This results has been improved
by De Lellis, Dominik Inauen and Szekelyhidi in 2018 where they proved
that any short immersion (or embedding) of a $2$-dimensional
Riemannian manifold $(B_1,g)$ into $\R^3$,
can be uniformly approximated by a sequence of $C^{1,\alpha}$ isometric
immersions (embeddings) for any $\alpha<\frac{1}{5}$. Their key
argument relies on the fact that every two-dimensional metric is
locally conformally equivalent to the Euclidean metric $g=\Id_2$.
This means that a positive definite defect $\mathcal D=g-(\nabla
u)^T\nabla u$ may be, by a change of variables, reduced to the diagonal
form, which decomposes into two primitive defects rather than three,
resulting in a lower rate of blow-up of the second derivatives in the Nash-Kuiper
iteration scheme and subsequently higher regularity of the immersions
derived in the limiting process. The same statement, albeit at the
linearized level, has been recently used by the second author of this
paper, to show density in the space of continuous functions on
$\bar\Omega\subset\R^2$, of the set of weak solutions $u\in
C^{1,\alpha}(\bar\Omega,\R^{n})$ with any $\alpha <\frac{1}{1+4/(n-2)}$,
to the following equation with a given right hand side $f$:
$$\partial_{11}u\cdot\partial_{22}u - |\partial_{12}u|^2=f \quad\text{in } \Omega.$$
For the codimension $n-2=1$, this generalizes the prior density
result for the {\em Monge-Amp\`ere equation} and its weak $C^{1,\alpha}$
solutions at $\alpha<\frac{1}{5}$, due to Wentao Cao and
Szekelyhidi. The parallel rigidity statements are likewise available
when $\alpha>\frac{2}{3}$.

\medskip

\noindent Regarding the flexibility vs rigidity in the regularity
interval $[\frac{1}{5},\frac{2}{3}]$, Gromov conjectured that the
actual threshold occurs sharply at
$\alpha=\frac{1}{2}$. This is supported by the work of De Lellis and Inauen
in 2020 in which they proved that for any $\alpha<\frac{1}{2}$,
an appropriate convex integration construction yields $C^{1,\alpha}$ isometric
immersions of a spherical cup, whose Levi-Civita connection differs from the
standard one, whereas any such immersion with regularity
$\alpha>\frac{1}{2}$ must necessarily induce the compatible Levi-Civita connection.

\section{Applications to general relativity}

\noindent {\em Quasi-local mass} in general relativity is a notion associated with
closed spacelike $2$-surfaces in a $4$-dimensional spacetime. Its
purpose is to evaluate the amount of matter and gravitational energy
contained within the surface, and can potentially be used to detect the formation of black holes.
In this section, we briefly discuss an application of
Weyl's embedding problem to the quasi-local masses. 
For more information, see \cite{Szab}.

\medskip

\noindent Consider a smooth, orientable, compact Riemannian
manifold $(\Omega^3,g)$, with connected boundary $\Sigma$ of positive Gaussian curvature.
According to Weyl's embedding theorem, $\Sigma$ may be uniquely (up
to rigid motions) isometrically embedded into $\mathbb R^3$. This
embedding induces the mean curvature $H_0$,
which, in general, differs from the mean curvature $H$ of $\Sigma$
as a submanifold of $\Omega$. In 1992, based on a Hamilton-Jacobi
analysis of the Einstein-Hilbert action, David Brown and James York
defined the quasi-local mass of $\Sigma$ to be: 
\begin{equation*}
m_{BY}(\Sigma) =\frac{1}{8\pi}\Big(\int_{\Sigma} H_0 \;\mbox{d}\sigma 
-\int_{\Sigma} H\;\mbox{d}\sigma\Big).
\end{equation*}
A fundamental result concerning $m_{BY}$ was established by Yuguang Shi and
Luen-Fai Tam in 2002. Namely, they showed that if 
$H$ is positive and the scalar curvature of $g$ is
nonnegative, then $m_{BY}(\Sigma)$ is nonnegative and
it vanishes if and only if $(\Omega,g)$ isometrically embeds into
$\mathbb{R}^3$. From a geometric perspective, this result may be
interpreted as a comparison theorem for compact manifolds of
nonnegative scalar curvature. 

\medskip

\noindent Despite this beautiful result, the Brown-York
definition has several deficiencies, most notably
that it is not `gauge independent' when considered in a spacetime
context. This motivated Chiu-Chu Melissa Liu and Shing-Tung Yau in
2003, and then Mu-Tao Wang and Yau in 2009 to each define more general notions of quasi-local
mass which satisfy a range of desirable properties.  
Like $m_{BY}$, both of these masses also employ the Weyl
embedding theorem, and are consequently restricted to surfaces
$\Sigma$ which are topologically 2-spheres. It should be noted that
Wang-Yau utilize the theorem to produce isometric embeddings into
Minkowski space, even if the Guassian curvature changes sign. Recently,
another quasi-local mass in this family was proposed by
Aghil Alaee, Khuri and Yau, which allows for surfaces of higher genus
and also requires their embedding into Minkowski space. A natural
question then arises that would have important implications: which
closed surfaces admit isometric embeddings
into Minkowski space? This problem gives rise to an underdetermined
system of equations, and thus one may guess that there are no
obstructions. However, even for the torus, the problem remains open.

\section{Applications to the mathematical materials science}

\noindent When the ambient and intrinsic dimensions agree, $n=d$, the problem
(\ref{A.1.4}) is linked with the satisfaction of the
orientation preservation by $u: B_1\to\R^d$, expressed as:
\begin{equation}\label{op_IN}
\det\nabla u>0 \quad \mbox{ in }\;B_1.
\end{equation}
Under this condition, a sufficient and necessary condition
for the local solvability of (\ref{A.1.4}) is the vanishing of the Riemannian curvature
of $g$, which also guarantees that the solution $u$ is smooth and
unique up to rigid motions. On the other hand, without the restriction (\ref{op_IN}), there always exists a Lipschitz
continuous $u$ constructed by convex integration, which indeed changes orientation in
any neighbourhood of any point at which $g$ has non-zero
curvature. The set of such Lipschitz immesions is dense in the set
of short immersions, similarly to other $h$-principle statements that we have listed before. 

\medskip

\noindent 
In the former context, it is natural to  pose the quantitative question: what is the infimum
of the averaged pointwise deficit of $u$ from being an
orientation-preserving isometric immersion of $g$? This deficit is
measured by the following {\em non-Euclidean energy} on a domain
$\Omega\subset\R^d$ with respect to the Riemannian manifold $(\Omega, g)$: 
\begin{equation}\label{ee_IN}
\mathcal{E}_g(u)=\int_\Omega \mathrm{dist}^2\big((\nabla u)g^{-1/2}, \mathrm{SO}(d)\big)\dx.
\end{equation}
Above, $\mathrm{SO}(d)$ denotes the special orthogonal group, namely rotations
in $\R^d$, and $\mathrm{dist}(\cdot, \cdot)$ is the distance in the space of matrices
$\R^{d\times d}$. By the polar decomposition theorem, an equidimensional
$u$ satisfies  (\ref{A.1.4}) and (\ref{op_IN}) if and only if 
$\nabla u \in \mathrm{SO}(d) g^{1/2}$ in $\Omega$,
which happens precisely when $\mathcal{E}_g(u)=0$. The follow-up questions now are:
\begin{enumerate}
\item[(i)] Can one quantify the infimum of $\mathcal{E}_g$ in relation to $g$ and $\Omega$? 
\item[(ii)] What is the structure of minimizers to (\ref{ee_IN}), if they exist? 
\item[(iii)] In the limit of $\Omega$ be\-co\-ming $(d-1)$-dimensional, what are
the asymptotic properties of the energies $\mathcal{E}_g$ and
their minimizers in relation to the curvatures of $g$?
\end{enumerate}

\subsection{Connection to calculus of variations and elasticity.}

The field of calculus of variations originally centered around minimization
problems for integral functionals of the general form below, in which
$W:\Omega\times \R^n\times\R^{n\times d}\to \R$ is a given
{\em energy density}, and where $u$ may be subject to various
constraints, for example the boundary conditions:
\begin{equation}\label{ge_IN}
\mathcal{E}(u) =  \int_\Omega W\big(x,u(x),\nabla u(x)\big)\dx \quad
\mbox{ for }\; u:\R^d\supset\Omega\to\R^n,
\end{equation}
The systematic study of existence of minimizers to (\ref{ge_IN}), their uniqueness and
qualitative properties, began with Leonhard Euler and Johann Bernoulli in the XVIIth century
and progressed due to seminal contributions by 
Charles Morrey and Ennio De Giorgi in the XXth century. These questions are strongly
tied to the convexity, in the $\nabla u$ variable in 
$W$, in turn implying the so-called sequential lower semicontinuity of $\mathcal{E}$;
a condition necessary to conclude that the minimizing sequences to
(\ref{ge_IN}) accumulate at the minimizers.
This is the {\em direct method of calculus of
  variations},  which however does not apply to the functional in
(\ref{ee_IN}), precisely due to its lack of convexity. 

\medskip

\noindent An example of an important class of problems of the form
(\ref{ge_IN}) is the basic variational
model pertaining to the {\em nonlinear elastic energy of prestressed bodies}:
\begin{equation}\label{none_IN}
\mathcal{E}_g(u) =  \int_\Omega W\big((\nabla u)g^{-1/2}\big)\dx \quad
\mbox{ for }\; u:\R^3\supset\Omega\to\R^3,
\end{equation}
where $W:\R^{3\times 3}\to\R$ is the given energy density, carrying the
elastic moduli of the physical material whose referential configuration
is $\Omega$, and satisfying the necessary physically-relevant conditions (frame
invariance and the zero-penalty for all rigid motions).
The theory of elasticity is one of the most important subfields of continuum
mechanics. It studies materials which are capable of undergoing large
deformations, due to the distribution of {\em local stresses and
displacements}, and resulting from the application of mechanical or thermal loads. 
The model (\ref{none_IN}) postulates formation of a
target Riemannian metric $g$ and the induced multiplicative decomposition of
the deformation gradient $\nabla u$ into the elastic part $(\nabla
u)g^{1/2}$, and the inelastic part $g^{1/2}$ responsible for the
morphogenesis. The form of $g$ is dictated by the material's response to pH, humidity, temperature,
growth hormone distribution and other stimuli, and it is specific to each problem.

\medskip

\noindent The functional in (\ref{none_IN}) corresponds to a range of
hyperelastic energies approxi\-ma\-ting the behavior of a
large class of elastomeric materials, and it is consistent with the microscopic
derivations based on statistical mechanics. It reduces (via a change of variables) to
the classical nonlinear three-dimensional elasticity, for $g$ with
vanishing Riemannian curvature, which occurs precisely when $\min \mathcal{E}_g=0$.
It can be proved that in the opposite case, i.e. for a non-Euclidean $g$, the infimum of
$\mathcal{E}_g$ in the absence of any
forces or boundary conditions remains strictly positive, pointing to the existence of residual strain. 

\medskip

\noindent For domains $\Omega$ that are {\em thin films}, one considers a family of problems:
\begin{equation}\label{eeh_IN}
\mathcal{E}^h_g(u^h)= \int_{\Omega^h} W\big((\nabla
u^h)g^{-1/2}\big)\dx \quad \mbox{ for }\; u:\Omega^h\to\R^3,
\end{equation}
parametrised by the small thickness $h$ of $\Omega^h=\omega\times (-h/2, h/2)$.
The task is now to determine the asymptotic limit of their minimizations
as $h\to 0$, rather than to minimize $\mathcal{E}_g^h$ for each
particular $h$. This can be achieved using the method of {\em
  $\Gamma$-convergence}, to identify the ``singular li\-mit'' energy
$\mathcal{I}_g$ characterized by the property that the minimizers and minimum 
values of (\ref{eeh_IN}) converge  to the minimizers and the minimum values
of $\mathcal{I}_g$. In general, one expects that $\inf
\mathcal{E}^h_g \simeq h^\gamma$ as $h\to 0$ with the 
optimal scaling exponent $\gamma$ determined from $g$ and the
appropriate curvatures of $g$ contributing
to the form of $\mathcal{I}_g$.

\subsection{Dimension reduction of thin prestressed films.}

Using flexibiliy of the isometric immersions of Riemannian manifolds $(M^2,g)$ into $\R^3$
at H\"older regularity below $C^{1,1/5}$, one can show that $\inf
\mathcal{E}_g^h\leq C h^\beta$ for any $\beta<\frac{2}{3}$. On the
other hand, having any $\beta\geq 2$ automatically implies (in fact,
it is equivalent to) that the restriction $g(\cdot, 0)_{2\times 2}$ of the metric $g$ to the
midplate $\omega\subset \R^2$ has an isometric immersion into $\R^3$ with Sobolev
regularity $H^2$ (i.e. with its second order derivatives square integrable).

\medskip

\begin{figure}[htbp]
\makebox[\textwidth][c]{
\renewcommand{\arraystretch}{1.5}
$ \begin{array}{| l | l | l | l |}
\hline \vspace{-7mm}\\ {\beta} & \mbox{\begin{minipage}{4cm} \footnotesize{constraint /
       regularity}\end{minipage}} &
\mbox{\begin{minipage}{3.2cm}  \footnotesize{asymptotic expansion} \end{minipage}} &
\begin{minipage}{2.8cm} \footnotesize{\mbox{limiting energy }
    ${\mathcal{I}}_{g}$ }\end{minipage} \\ [10pt]\hline
\hspace{0.5mm}{2} & \hspace{-1.5mm} \footnotesize{ \begin{array}{l} y\in H^{2}(\omega,\R^3)
\\ (\nabla y)^T\nabla y = g(\cdot,0)_{2\times 2} \end{array}} 
& \hspace{-1.5mm} \footnotesize{\begin{array}{l} y(\cdot) \\ \hspace{-1mm} 
\big\{ 3d:~ y + x_3 b \big\}  \end{array} }
& \hspace{-2mm} \footnotesize{\begin{array}{l} \|(\nabla
y)^T\nabla b - \frac{1}{2}\partial_3g(\cdot,0)_{2\times
  2}\|^2 \vspace{0.5mm} \\ \big[\partial_1y,\partial_2y, b\big]\in
\mathrm{SO}(3)g(\cdot,0)^{1/2}\end{array} }\\ [15pt]\hline 
\hspace{0.5mm}{4} & \hspace{-2mm} \footnotesize{\begin{array}{l}
  \mathcal{R}_{12,cd}(\cdot,0) =0 \hspace{-2mm} 
\\ (V, w^h)\in H^{2}\times H^1(\omega, \R^3)
\\ \sym\big((\nabla y_0)^T\nabla 
  V\big)=0, \\  \sym \big((\nabla y_0)^T\nabla w^h\big)\to {S} \end{array} } \hspace{-2mm} 
& \hspace{-1.5mm} \footnotesize{ \begin{array}{l} y_0+hV +\; h^2w^h\end{array} }
& \hspace{-2mm} \footnotesize{\begin{array}{l}
  \|\frac{1}{2}(\nabla V)^T\nabla V +  {S} + \frac{1}{24}(\nabla b_1)^T\nabla  b_1
  \\ \hspace{4mm} -\frac{1}{48}\partial_{33}g(\cdot,0)_{2\times 2}\|^2 \\ +
\|(\nabla y_0)^T\nabla p + (\nabla V)^T\nabla b_1\|^2 \\  
+ \|\big[\mathcal{R}_{i3,j3}(\cdot,0)\big]_{i,j=1,2}\|^2 \vspace{2mm} \end{array}}
\\ [15pt] \hline  
\begin{array}{c} \hspace{-1mm} {6} \\ \hspace{-2mm} {\vdots} \end{array} 
&\hspace{-2mm} \footnotesize{\begin{array}{l} \mathcal{R}_{ab,cd}(\cdot,0)=0\\ 
V\in H^{2}(\omega,\R^3) \\ \sym \big((\nabla y_0)^T\nabla V\big)=0\hspace{-2mm}
\vspace{1mm}\end{array} } 
& \hspace{-1.5mm} \footnotesize{ \begin{array}{l} y_0 + h^2V \end{array}}
& \hspace{-2mm} \footnotesize{\begin{array}{l}
\|(\nabla y_0)^T\nabla p + (\nabla V)^T\nabla  
  b_1 +\alpha \big[\partial_3 \mathcal{R}\big] \|^2
\\ + \|\mathbb{P}_{\mathcal{S}_{y_0}^\perp}\big[\partial_3\mathcal{R}\big]\|^2
+ \|\mathbb{P}_{\mathcal{S}_{y_0}}\big[\partial_3 \mathcal{R}\big]\|^2 
\end{array} } \\ [15pt]\hline 
\hspace{-1.5mm}\begin{array}{c}  {2k} \\
 \hspace{-2mm}  {\vdots} \end{array} \hspace{-1.5mm} 
& \hspace{-2mm} \footnotesize{ \begin{array}{l} \mathcal{R}_{ab,cd}(\cdot,0)=0 \\ 
  \big[\partial_3^{(i)}\mathcal{R}\big](\cdot,0) = 0 \;\;\forall i\leq n-3 \hspace{-2mm} \\ 
V\in H^{2}(\omega,\R^3) \\ \sym \big((\nabla y_0)^T\nabla
V\big)=0 \vspace{1mm}\end{array} \hspace{-2mm} }
& \hspace{-1.5mm} \footnotesize{\begin{array}{l} y_0+h^{k-1}V \\ 
\hspace{-1mm}\big\{ 3d:~ y_0 + \displaystyle{\sum_{i=1}^{k-1} }\frac{x_3^i}{i!}
  b_i  \\ \hspace{4mm} + h^{k-1}V +
  h^{k-1} x_3 p\big\}  \end{array}} 
& \hspace{-2mm} \footnotesize{ \begin{array}{l}
 \|(\nabla y_0)^T\nabla p + (\nabla V)^T\nabla  b_1
\\ \hspace{4mm} + \alpha \big[\partial_3^{(k-2)} \mathcal{R}\big] \|^2 
+ \|\mathbb{P}_{\mathcal{S}_{y_0}^\perp}\big[\partial_3^{(k-2)}\mathcal{R}\big]\|^2
\hspace{-4mm} \\
+ \|\mathbb{P}_{\mathcal{S}_{y_0}}\big[\partial_3^{(k-2)}R\big]\|^2 
\end{array} } \\ [10pt]\hline \end{array} $ }
\caption{The first column gathers equivalent conditions for the
  scaling $\inf \mathcal{E}_g^h\sim Ch^\beta$, in terms of the Riemann
  curvatures $\mathcal{R}_{ab,cd}$ of $g$. The second column provides
  the asymptotic expansion of the minimizing sequence. The third
  column gives the form of the $\Gamma$-limits in this infinite hierarchy.}
\label{table_prestrain}
\end{figure}

\noindent One can show that the only viable scaling exponents in: $\inf
\mathcal{E}_g^h\sim Ch^\beta$ in the regime $\beta\geq 2$, are the even powers $\beta = 2k$.
For each $k$, the complete set of results is available:
the conditions/constraints on the curvatures of $g$ equivalent to the
indicated scaling; the necessary asymptotic expansion of the minimizing sequence
$u^h$ in terms of the transversal (thin direction) variable $x_3$ and
the regularity of the fields present in that expansion; the limiting
energy $\mathcal{I}_g$ given in terms of these fields and the
unconstrained curvatures. That {\em hierarchy of dimensionally reduced energies}
is schematically presented in Figure \ref{table_prestrain}.

\medskip

\noindent Parallel general results can be derived in the abstract setting of
Riemannian manifolds, for more general dimensions of the midplate
$d>2$ and the ambient space $n>d$, and also for $g$ replaced by the film's
thickness - depending prestrain, incompatible only through a
perturbation of order which is a power of $h$. 
While the systematic description of the singular limits associated with the
exponents $\beta\in (\frac{2}{3},2)$ is not yet available, there are a
number of illustrative examples of the emerging patterns and the
corresponding scalings in that range. For this discussion we
refer to \cite{lew_book} and the references therein.

\end{document}